\documentclass[12pt]{amsart}
\usepackage{amssymb}
\usepackage{amsfonts}
\usepackage{color,soul}
\usepackage{graphicx}
\usepackage{relsize}

\newcommand{\nc}{\newcommand}
\nc{\thref}[1]{Theorem~\ref{theo:#1}}
\nc{\selabel}[1]{\label{sect:#1}}
\nc{\seref}[1]{Section~\ref{sect:#1}}
\nc{\lelabel}[1]{\label{lemm:#1}}
\nc{\leref}[1]{Lemma~\ref{lemm:#1}}
\nc{\prlabel}[1]{\label{prop:#1}}
\nc{\prref}[1]{Proposition~\ref{prop:#1}}
\nc{\colabel}[1]{\label{coro:#1}}
\nc{\coref}[1]{Corollary~\ref{coro:#1}}
\nc{\exlabel}[1]{\label{exam:#1}}
\nc{\exref}[1]{Example~\ref{exam:#1}}
\nc{\delabel}[1]{\label{defi:#1}}
\nc{\deref}[1]{Definition~\ref{defi:#1}}
\nc{\eqlabel}[1]{\label{equa:#1}}
\nc{\relabel}[1]{\label{rema:#1}}
\nc{\reref}[1]{Lemma~\ref{rema:#1}}
\providecommand{\operatorname}[1]{\mathrm{#1}\,}
\nc{\Hom}{\operatorname{Hom}} \nc{\Mor}{\operatorname{Mor}}
\nc{\Aut}{\operatorname{Aut}} \nc{\Ann}{\operatorname{Ann}}
\nc{\Ker}{\operatorname{Ker}} \nc{\Trace}{\operatorname{Trace}}
\nc{\Char}{\operatorname{Char}} \nc{\Mod}{\operatorname{Mod}}
\nc{\End}{\operatorname{End}} \nc{\Spec}{\operatorname{Spec}}
\nc{\Span}{\operatorname{Span}} \nc{\sgn}{\operatorname{sgn}}
\nc{\Id}{\operatorname{Id}} \nc{\Com}{\operatorname{Com}}
\nc{\rank}{\operatorname{rank}}
\nc{\Clausen}{\operatorname{Cl}}

\let\:=\colon

\newtheorem{de}{Definition}[section]
\newtheorem{lm}[de]{Lemma}
\newtheorem{pr}[de]{Proposition}
\newtheorem{co}[de]{Corollary}
\newtheorem{re}[de]{Remark}
\newtheorem{res}[de]{Remarks}
\newtheorem{te}[de]{Theorem}
\newtheorem{ex}[de]{Example}
\newtheorem{exs}[de]{Examples}

\def\bex{\begin{ex}}
\def\eex{\end{ex}}
\def\bexs{\begin{exs}}
\def\eexs{\end{exs}}
\def\bl{\begin{lm}}
\def\el{\end{lm}}
\def\bc{\begin{co}}
\def\ec{\end{co}}
\def\bt{\begin{te}}
\def\et{\end{te}}
\def\bpr{\begin{pr}}
\def\epr{\end{pr}}
\def\br{\begin{re}}
\def\er{\end{re}}
\def\brs{\begin{res}}
\def\ers{\end{res}}
\def\bd{\begin{de}}
\def\ed{\end{de}}
\def\be{\begin{equation}}
\def\ee{\end{equation}}
\def\bea{\begin{eqnarray*}}
\def\eea{\end{eqnarray*}}
\def\bp{\begin{proof}}
\def\ep{\end{proof}}

\def\qed{\hfill\Box}

\def\RR{{\mathbb R}}

\def\NN{{\mathbb N}}

\nc{\starint}{\mathlarger{\ast}\kern-0.83em{\int}}

\let\:=\colon

\begin{document}

\title{Redefining the integral}

\author{Derek Orr}

\begin{abstract}
In this paper, we discuss a similar functional to that of a standard integral. The main difference is in its definition: instead of taking a sum, we are taking a product. It turns out this new ``star-integral" may be written in terms of the standard integral but it has many different (and similar) interesting properties compared to the regular integral. Further, we define a ``star-derivative" and discuss its relationship to the ``star-integral".
\end{abstract}



\maketitle

\section{Introduction}

The definition of an integral has been around since the 17th century and even earlier when mathematicians such as Cavalieri and Fermat began experimenting with integrals of polynomials (see \cite{Ander}). The definite integral has the precise meaning of the area underneath a graph $y=f(x)$. 

Formally, given an interval $[a,b]$, there exists a partition of it $P_{n}(a,b) := \{[x_{i-1},x_{i}] : a = x_{0} < x_{1} < x_{2} < \ldots < x_{n-1} < x_{n} = b\}$ with the length of the intervals $\Delta x_{i} := x_{i}-x_{i-1}$, and a test set $T(P):= \{t_{i} : x_{i-1} < t_{i} < x_{i}, \hspace{3pt} \forall \hspace{3pt} i\in\{1,2,\ldots,n\}\}$ where $\big(t_{i}\big)_{i\in\{1,2,\ldots,n\}}$ form a sequence of ``test points". Further if we define $\displaystyle \mathcal{S}(f,P,T) := \sum_{i=1}^{n} f(t_{i})\Delta x_{i}$ and the mesh of the partition, $\displaystyle mesh(P):=\max_{i\in\{1,\ldots,n\}} \Delta x_{i}$, we can state the definition of a function being Riemann integrable.

\vspace{0.5cm}

\bd[Riemann integral]{For $S\in\RR$, if $\forall \hspace{3pt} \epsilon > 0, \hspace{3pt} \exists \hspace{3pt} \delta > 0$  such that}

$$mesh(P) < \delta \hspace{4pt} \Longrightarrow \hspace{4pt} \Bigg|\mathcal{S}(f,P,T) - S\Bigg| < \epsilon, $$

then $f$ is Riemann integrable with Riemann integral $S$ and write $\displaystyle \int_{a}^{b} f(x) \hspace{3pt} dx = S$.

\ed

\vspace{0.5cm}
Now, let's consider an arbitrary Riemann integrable function. Usually the partition will be an equipartition, which means $\Delta x_{i} = \Delta x = (b-a)/n, \hspace{3pt} \forall \hspace{3pt} i\in\{1,2,\ldots,n\}$ and then the definition above becomes a bit simpler. Further, for simplicity, let $t_{i} = x_{i} = a+i\Delta x$ here. Thus, we will have $\displaystyle \mathcal{S}(f,P,T) = \sum_{i=1}^{n} f(a+i\Delta x)\Delta x$ and since $f(x)$ is Riemann integrable,

$$\lim_{n\to\infty} \mathcal{S}(f,P,T) = \lim_{n\to\infty} \frac{b-a}{n}\sum_{i=1}^{n}f\Big(a+i\frac{b-a}{n}\Big) = S.$$

For the purposes of this paper, we will use this as the method to compute an integral. 

\vspace{0.5cm}

\section{The star-integral}

Partitioning an interval into equal-sized subintervals and taking the length of the subintervals to be small has been an incredibly common technique to figure out the value of an integral. But if we think about this more abstractly, we can define a new type of integral. The main difference will be instead of a sum, we will have a product. Infinite products are far less common than infinite sums but they still appear in many areas of mathematics such as number theory involving the Riemann zeta function or gamma function (see \cite{Borwein}--\cite{Ram}). So let us consider

$$\lim_{n\to\infty} \frac{b-a}{n}\prod_{i=1}^{n} f\Big(a+i\frac{b-a}{n}\Big).$$

\vspace{0.2cm}

A good first question is if this is well-defined. For a fixed $n$, this product is finite, since $f(x)$ is assumed to be well-behaved on $[a,b]$. But, some things are definitely not right here. For example, the first term in the product is $\displaystyle f\Big(\frac{b}{n}+a\Big(1-\frac{1}{n}\Big)\Big)$. So as $n$ gets larger, the first term approaches $f(a)$. If $f(a)=0$, this infinite product could approach 0 which could ruin a lot of analysis for this expression. We must change more than just the summation sign. Keep in mind, a product is essentially an iterated sum. For integers, one can think of $nx = x+x+\ldots+x$ a total of $n$ times. So, by switching the sum to a product, we've switched the sum to an iterated sum. So, with $f(x)$ and $dx$ multiplying each other in $\mathcal{S}$, we need to switch that operation to an ``iterative multiplication" operation. The clear choice here is exponentiation. So we will denote $\displaystyle \mathcal{P}(f,P,T) := \prod_{i=1}^{n} f(t_{i})^{\Delta x_{i}}$ and we will have the following:

\bd[Riemann star-integral]{For $M\in\RR$, if $\forall \hspace{3pt} \epsilon > 0, \hspace{3pt} \exists \hspace{3pt} \delta > 0$  such that}

$$mesh(P) < \delta \hspace{4pt} \Longrightarrow \hspace{4pt} \Bigg|\mathcal{P}(f,P,T) - M\Bigg| < \epsilon, $$

we say $f$ is Riemann star-integrable with Riemann star-integral $M$ and we write $\displaystyle \starint_{\hspace{-5pt}a}^{b} f(x) \hspace{3pt} dx = M$.

\ed

\vspace{0.5cm}

And for the equipartition, we have

$$\lim_{n\to\infty} \mathcal{P}(f,P,T) = \lim_{n\to\infty} \Bigg(\prod_{i=1}^{n}f\Big(a+i\frac{b-a}{n}\Big)\Bigg)^{\frac{b-a}{n}} = M.$$

\vspace{0.3cm}

Let's try to use this in the following example.

\bex{} Compute

$$\starint_{\hspace{-5pt}0}^{1} x \hspace{3pt} dx. $$

\eex

\textit{Solution.} Using equipartitions, we will have

$$\starint_{\hspace{-5pt}0}^{1} x \hspace{3pt} dx = \lim_{n\to\infty}\Bigg(\prod_{i=1}^{n}\Big(\frac{i}{n}\Big)\Bigg)^{1/n} = \lim_{n\to\infty} \bigg(\frac{n!}{n^n}\bigg)^{1/n}$$

\vspace{0.3cm}

and using Stirling's formula $\displaystyle \lim_{n\to\infty} \frac{n!e^n}{n^n\sqrt{2\pi n}} = 1$,

$$\starint_{\hspace{-5pt}0}^{1} x \hspace{3pt} dx = \lim_{n\to\infty} \bigg(\frac{\sqrt{2\pi n}}{e^n}\bigg)^{1/n} = \frac{1}{e}.$$

\vspace{0.3cm}

In fact, one can show $\displaystyle \starint_{\hspace{-5pt}0}^{x} t \hspace{3pt} dt = \Big(\frac{x}{e}\Big)^x$, which is very different from $\displaystyle \int_{0}^{x} t \hspace{3pt} dt = \frac{1}{2}x^2$. Although there are some integrals which make sense. For example, it's not hard to see that $\displaystyle \starint_{\hspace{-5pt}a}^{a} f(x) \hspace{3pt} dx = 1$. This is intuitive because the regular integral of a function with $b=a$ is 0, the additive identity, whereas the star-integral of a function with $b=a$ is 1, the multiplicative identity. Also, $\displaystyle \starint_{\hspace{-5pt}a}^{b} 1 \hspace{3pt} dx = 1$ is analogous to $\displaystyle \int_{a}^{b} 0 \hspace{3pt} dx = 0$.

\vspace{0.5cm}

\section{Connection to standard integral}

This is pretty interesting but we need to see how it acts on other functions. Also, it is time consuming to keep using the product definition of a star-integral. Is there a way we can write this star-integral in terms of the regular integral we know a lot about? If we go back to the definition of $\mathcal{P}(f,P,T)$, the main way to relate a product to a sum is by using logarithms, namely the natural logarithm which we will denote by $\log(x)$. The main property we will use is $\forall \hspace{3pt} x,y > 0, \log(xy) = \log(x) + \log(y)$. This can then be iterated to show the logarithm of a product is the sum of the logarithms for a finite number or countably infinite number of terms. It is around this time here that we mention we will be taking $f(x) > 0$ on the domain of star-integration and $f(x)=0$ only on sets of measure zero. Now, let's apply the logarithm to $\mathcal{P}$ to find

$$\log\big(\mathcal{P}(f,P,T)\big) = \log\Big(\prod_{i=1}^{n} f(t_{i})^{\Delta x_{i}}\Big) = \sum_{i=1}^{n} \Delta x_{i} \log(f(t_{i})) = \mathcal{S}(\log(f(x)),P,T),$$ 

\vspace{0.2cm}

where we have also used the fact that $\log(x^y) = y\log(x)$ for $x>0$ and $y\in \RR$. Since $\log(x)$ is a continuous function, taking limits gives

\begin{equation}
    \starint_{\hspace{-5pt}a}^{b} f(x) \hspace{3pt} dx = \exp{\bigg(\int_{a}^{b} \log(f(x)) \hspace{3pt} dx\bigg)},
\end{equation}

\vspace{0.2cm}

and now we see an obvious candidate for the indefinite star-integral:

\begin{equation}
    \starint f(x) \hspace{3pt} dx = \exp{\bigg(\int \log(f(x)) \hspace{3pt} dx\bigg)}.
\end{equation}

\vspace{0.3cm}

Now we can examine the properties of star-integrals and how they are related to the properties of regular integrals. Here are a few:

\brs{Properties of the star-integral, provided all integrals mentioned exist:}

\vspace{0.2cm}

(a) $\hspace{10pt} \displaystyle f(x) \leq g(x) \implies \starint_{\hspace{-3pt}a}^{b} f(x) \hspace{3pt} dx \leq \starint_{\hspace{-3pt}a}^{b} g(x) \hspace{3pt} dx$

\vspace{0.2cm}

(b) $\hspace{10pt} \displaystyle \starint e^{f(x)} \hspace{3pt} dx = \exp\Bigg(\int f(x) \hspace{3pt} dx\Bigg) \Longleftrightarrow \log\Bigg(\starint f(x) \hspace{3pt} dx\Bigg) = \int \log(f(x)) \hspace{3pt} dx$

\vspace{0.2cm}

(c) $\hspace{10pt} \displaystyle \starint f(x)^ng(x)^m \hspace{3pt} dx = \Bigg(\starint f(x) \hspace{3pt} dx\Bigg)^n\Bigg(\starint g(x) \hspace{3pt} dx\Bigg)^m, \hspace{5pt} m,n \in \RR$

\vspace{0.2cm}

(d) $\hspace{10pt} \displaystyle \starint c f(x) \hspace{3pt} dx = c^x\hspace{2pt}\starint f(x) \hspace{3pt} dx, \hspace{5pt} c > 0$

\vspace{0.2cm}

(e) $\hspace{10pt} \displaystyle \starint_{\hspace{-5pt}b}^{a} f(x) \hspace{3pt} dx = \starint_{\hspace{-5pt}a}^{b} f(x)^{-1} \hspace{3pt} dx$

\vspace{0.2cm}

(f) $\hspace{10pt} \displaystyle \starint_{\hspace{-5pt}a}^{c} f(x) \hspace{3pt} dx = \Bigg(\starint_{\hspace{-5pt}a}^{b} f(x) \hspace{3pt} dx\Bigg)\Bigg(\starint_{\hspace{-5pt}b}^{c} f(x) \hspace{3pt} dx\Bigg)$

\vspace{0.2cm}

(g) $\hspace{10pt} \displaystyle \starint_{\hspace{-5pt}a}^{b} s f(x)+ (1-s) g(x) \hspace{3pt} dx \geq \Bigg(\starint_{\hspace{-5pt}a}^{b} f(x)^s \hspace{3pt} dx\Bigg)\Bigg(\starint_{\hspace{-5pt}a}^{b} g(x)^{1-s} \hspace{3pt} dx\Bigg), \hspace{5pt} s \in [0,1]$

\ers

\vspace{0.2cm}

The first property holds because $\exp(x)$ and $\log(x)$ are increasing functions. The second property is entirely clear by $(2)$. Formulas 3.1.c and 3.1.d use additivity of the regular integral and properties of $\log(x)$ and $\exp(x)$. The next two are straightforward consequences of $(1)$. The last one is quite interesting, almost a sort of H\"older inequality. Clearly by $(c)$, standard H\"older's inequality is gone. But by the concavity of $\log(x)$ and the fact that $\exp(x)$ is increasing, we arrive at $(g)$, where equality occurs if and only if $f(x)=g(x)$. And for $s=1/2$, using the above remarks, $(g)$ can be transformed to

\begin{equation}
\starint_{\hspace{-5pt}a}^{b} f(x) + g(x) \hspace{3pt} dx \geq 2^{b-a}\Bigg(\starint_{\hspace{-5pt}a}^{b} f(x)g(x) \hspace{3pt} dx\Bigg)^{1/2}
\end{equation}

\vspace{0.3cm}

which is perhaps the Cauchy-Schwarz inequality for star-integrals. Notice when $f(x)$ and $g(x)$ are constants, $(3)$ turns into the AM-GM inequality for 2 variables. 

\bl{$\forall \hspace{3pt} \alpha,\beta,\gamma > 0$,}

\begin{equation}
    (\alpha+\beta+\gamma)^2 \geq 8\alpha\sqrt{\beta\gamma}.
\end{equation}

\el

\vspace{0.3cm}

\textit{Proof.} Beginning with

$$\starint_{\hspace{-5pt}a}^{b} f(x)+g(x)+h(x) \hspace{3pt} dx \geq 2^{b-a}\Bigg(\starint_{\hspace{-5pt}a}^{b}f(x)\big(g(x)+h(x)\big) \hspace{3pt} dx\Bigg)^{1/2} = 2^{b-a}\Bigg(\starint_{\hspace{-5pt}a}^{b}f(x)g(x)+f(x)h(x) \hspace{3pt} dx\Bigg)^{1/2}.$$

\vspace{0.2cm}

Applying $(3)$ again, we have

\begin{multline*}
\starint_{\hspace{-5pt}a}^{b} f(x)+g(x)+h(x) \hspace{3pt} dx \geq 2^{\frac{3(b-a)}{2}}\Bigg(\starint_{\hspace{-5pt}a}^{b}f(x)g(x)f(x)h(x) \hspace{3pt} dx\Bigg)^{1/4} \\= 2^{\frac{3(b-a)}{2}}\Bigg(\starint_{\hspace{-5pt}a}^{b}f(x)^{1/2}\hspace{3pt} dx\Bigg) \Bigg(\starint_{\hspace{-5pt}a}^{b} g(x)^{1/4} \hspace{3pt} dx\Bigg)\Bigg(\starint_{\hspace{-5pt}a}^{b} h(x)^{1/4} \hspace{3pt} dx\Bigg).
\end{multline*}

\vspace{0.3cm}

Letting $f(x)=\alpha$, $g(x) = \beta$, $h(x) = \gamma$, and using $(1)$, we see

$$(\alpha+\beta+\gamma)^{b-a} \geq 2^{\frac{3(b-a)}{2}}\Big(\alpha^{\frac{b-a}{2}}\Big)\Big(\beta^{\frac{b-a}{4}}\Big)\Big(\gamma^{\frac{b-a}{4}}\Big).$$

\vspace{0.2cm}

Cancelling the $b-a$ in each exponent and squaring both sides, we arrive at the result. $\qed$

\vspace{0.4cm}

For best results, one should take $\alpha \geq \beta \geq \gamma$. In fact, this lemma can be generalized to the following theorem.

\bt{$\forall \hspace{3pt} a_{i} >0$ such that $a_{1} \leq a_{2} \leq a_{3} \leq \ldots \leq a_{k+1}$ and $\forall \hspace{3pt} k \in \NN$,}

\begin{equation}
    \Bigg(\sum_{i=1}^{k+1} a_{i}\Bigg)^{2^{k-1}} \geq 2^{2^k-1}\sqrt{a_{1}}\prod_{i=1}^{k}a_{i+1}^{2^{i-2}}.
\end{equation}

\et

\vspace{0.4cm}

The idea behind the proof is to use $(3)$ multiple times. Of course, $(5)$ still holds true in the case when $k=0$ or if some $a_{i}=0$. It is interesting to note that equality does not necessarily occur when all $a_{i}$ are equal. For example, setting $a_{i}=1$, we have

$$(k+1)^{2^{k-1}} \geq 2^{2^k-1}, \hspace{6pt} \forall \hspace{4pt} k \in \NN.$$

\vspace{0.2cm}

This even holds true $\forall \hspace{3pt} k \geq 1$. Going further, one can use the concavity of $\log(x)$ and generalize $(g)$ to

$$\starint_{\hspace{-5pt}a}^{b}\sum_{i=1}^{n}t_{i}f_{i}(x) \hspace{3pt} dx \geq \prod_{i=1}^{n}\Bigg(\starint_{\hspace{-5pt}a}^{b}f_{i}(x)^{t_{i}}\hspace{3pt} dx\Bigg), \hspace{7pt} t_{i} \geq 0, \hspace{4pt} \sum_{i=1}^{n} t_{i} = 1.$$

\vspace{0.3cm}

By setting all $t_{i} = 1/n$, we have

$$\starint_{\hspace{-5pt}a}^{b}\sum_{i=1}^{n}f_{i}(x) \hspace{3pt} dx \geq n^{b-a}\Bigg(\starint_{\hspace{-5pt}a}^{b}\prod_{i=1}^{n}f_{i}(x)\hspace{3pt} dx\Bigg)^{1/n}$$

\vspace{0.3cm}

and if one sets $f_{i}(x)=a_{i}>0$, we arrive at the AM-GM inequality for $n$ variables. 

\vspace{0.4cm}

\section{Star-integral formulas}

\vspace{0.3cm}

Now let's start to figure out some explicit formulas for the star-integral. What is the star-integral of $x^n$ or $e^x$ or other functions?

\bc{We have the following formulas, with $C$ an arbitrary constant:} 

\vspace{0.2cm}

(a) $\displaystyle \starint x^n \hspace{3pt} dx = C\bigg(\frac{x}{e}\bigg)^{nx}, \hspace{5pt} n \in \RR$

\vspace{0.2cm}

(b) $\displaystyle \starint e^x \hspace{3pt} dx = Ce^{x^2/2}$

\vspace{0.2cm}

(c) $\displaystyle \starint e^{k/x} \hspace{3pt} dx = Cx^k, \hspace{5pt} k \in \RR$

\vspace{0.2cm}

(d) $\displaystyle \starint e^{e^x} \hspace{3pt} dx = Ce^{e^x}$

\vspace{0.2cm}

(e) $\displaystyle \starint x^x \hspace{3pt} dx = C\bigg(\frac{x^2}{e}\bigg)^{x^2/4}$

\vspace{0.2cm}

(f) $\displaystyle \starint a^x \hspace{3pt} dx = Ca^{x^2/2},\hspace{5pt} a >0$

\vspace{0.2cm}

(g) $\displaystyle \starint ax+b \hspace{3pt} dx = C\bigg(\frac{ax+b}{e}\bigg)^x(ax+b)^{b/a}, \hspace{5pt} a,b \in \RR$

\ec

\vspace{0.3cm}

Since we have an indefinite integral in the exponent of $e^x$, the star-integral's arbitrary constant comes out as a coefficient of the function. All of these formulas are not hard to see by using $(2)$. What is interesting is we are able to easily compute star-integrals of difficult functions such as $(c)$, $(d)$, and $(e)$, but computing star-integrals of easier functions, like $(g)$, can be challenging. 

This also makes one wonder about the set of all Riemann star-integrable functions with a closed expression. Let $\mathcal{R}_{c}$ denote the standard set of Riemann integrable functions with a closed expression and let $\mathcal{R}_{c}^{*}$ denote the set of Riemann star-integrable functions with a closed expression. It's clear that they have the relation of $f(x) \in \mathcal{R}_{c} \Leftrightarrow e^{f(x)} \in \mathcal{R}_{c}^{*}$ and similarly $f(x) \in \mathcal{R}_{c}^{*} \Leftrightarrow \log(f(x)) \in \mathcal{R}_{c}$. But, neither is a subset of the other. For example, $e^{x^2} \in \mathcal{R}_{c}^{*}$ but it is famously known that $e^{x^2} \notin \mathcal{R}_{c}$. And $\log(x) \in \mathcal{R}_{c}$ but not in $\mathcal{R}_{c}^{*}$. Another question could be do all functions lie in $\mathcal{R}_{c} \cup \mathcal{R}_{c}^{*}$? It seems not; $x^{x^x}$ does not appear to be in either set. We could then ignore closed expressions and investigate simply a function being Riemann star-integrable compared to being Riemann integrable.

\section{Star-derivative}

\vspace{0.3cm}

There are still a few important rules to cover, namely the fundamental theorem of calculus. Here, the regular derivative does not help because $\displaystyle \starint f'(x) \hspace{3pt} dx$ does not simplify in any nice way. So to fix this, we need a different derivative. Also, is there an integration by parts formula? If we use the regular derivative function, we will find

$$\starint f(x)^{g'(x)} \hspace{3pt} dx = f(x)^{g(x)}\starint e^{-g(x)f'(x)/f(x)} \hspace{3pt} dx $$

\vspace{0.3cm}

but maybe there is a better formula? And again, $\displaystyle \starint f(x)g'(x) \hspace{3pt} dx$ does not simplify by formula 3.1.c. Recall the definition of a derivative below.

\vspace{0.4cm}

\bd[Definition of a derivative]{If $f$ is differentiable at a point $x$, we denote this as $f'(x)$ or $\displaystyle \frac{df}{dx}$ and}

$$\frac{df}{dx} = f'(x) = \lim_{h\to0} \frac{f(x+h)-f(x)}{h}. $$

\ed

\vspace{0.3cm}

Seeing how we changed the definition of the integral, we will do the same here.

\vspace{0.3cm}

\bd[Definition of a star-derivative]{If $f$ is star-differentiable at a point $x$, we denote this as $f^{*}(x)$ or $\displaystyle \frac{\delta f}{\delta x}$ and}

$$\frac{\delta f}{\delta x} = f^{*}(x) = \lim_{h\to0} \Bigg(\frac{f(x+h)}{f(x)}\Bigg)^{1/h} $$

\ed{}

\vspace{0.3cm}

where we've changed subtraction to division and multiplication to exponentiation. Before we unpack this, let us compute a couple of star-derivatives. It makes sense that $\displaystyle \frac{\delta}{\delta x}\big(C\big) = 1$, the multiplicative identity because the regular derivative of a constant is 0, the additive identity. This also agrees with what we noticed about star-integrals. Also 

$$\frac{\delta}{\delta x}\big(x\big) = \lim_{h\to0} \Bigg(\frac{x+h}{x}\Bigg)^{1/h} = \lim_{n\to\infty} \Bigg(1+\frac{1/x}{n}\Bigg)^n = e^{1/x}.$$

\vspace{0.3cm}

Now let's unpack this formula and take the natural logarithm of both sides. By continuity, we have

$$\log(f^{*}(x)) = \lim_{h\to0} \log\Bigg(\frac{f(x+h)}{f(x)}\Bigg)^{1/h} = \lim_{h\to0} \frac{\log(f(x+h)-\log(f(x))}{h} = \frac{d}{dx}\Big(\log(f(x))\Big)$$

\vspace{0.4cm}

and so

\begin{equation}
    f^{*}(x) = \exp\Bigg(\frac{d}{dx}\Big(\log(f(x))\Big)\Bigg) = \exp\Bigg(\frac{f'(x)}{f(x)}\Bigg).
\end{equation}

\vspace{0.4cm}

Notice how similar this is to the star-integral. To show the similarities, we give remarks to those in section 3.

\brs{Properties of the star-derivative}

\vspace{0.2cm}

(a) $\hspace{10pt} \displaystyle \frac{\delta}{\delta x}\Big(e^{f(x)}\Big) = \exp(f'(x)) \Longleftrightarrow \log(f^{*}(x)) = \frac{d}{dx}\Big(\log(f(x))\Big)$

\vspace{0.2cm}

(b) $\hspace{10pt} \displaystyle \frac{\delta}{\delta x}\Big(f(x)^ng(x)^m\Big) = \Big(f^{*}(x)\Big)^n\Big(g^{*}(x)\Big)^m, \hspace{5pt} m,n \in \RR$

\vspace{0.2cm}

(c) $\hspace{10pt} \displaystyle \frac{\delta}{\delta x} \Big(c f(x)\Big) = f^{*}(x), \hspace{5pt} c > 0$

\vspace{0.2cm}

(d) $\hspace{10pt} \displaystyle \frac{\delta}{\delta x}\Big(f(x)^{g(x)}\Big) = f(x)^{g'(x)}f^{*}(x)^{g(x)}$

\ers

\vspace{0.4cm}

Formula (d) is the product rule (or here, the ``exponent rule") for star-integrals, which hints at the integration by parts formula coming up. Moving on, we also compute some explicit star-derivatives below.

\bc{We have the following formulas} 

\vspace{0.2cm}

(a) $\displaystyle \frac{\delta}{\delta x}\Big(x^n\Big) = e^{n/x}, \hspace{5pt} n \in \RR$

\vspace{0.2cm}

(b) $\displaystyle \frac{\delta}{\delta x}\Big(e^x\Big) = e$

\vspace{0.2cm}

(c) $\displaystyle \frac{\delta}{\delta x} \Big(e^{e^x}\Big) = e^{e^x}$

\vspace{0.2cm}

(d) $\displaystyle \frac{\delta}{\delta x}\Big(x^x\Big) = ex $

\vspace{0.2cm}

(e) $\displaystyle \frac{\delta}{\delta x} \Big((nx)^x\Big)  = enx, \hspace{5pt} n > 0$

\vspace{0.2cm}

(f) $\displaystyle \frac{\delta}{\delta x} \Big(a^x\Big) = a,\hspace{5pt} a >0$

\vspace{0.2cm}

(g) $\displaystyle \frac{\delta}{\delta x}\Big( ax+b\Big) = e^{a/(ax+b)}, \hspace{5pt} a,b \in \RR$

\ec

\vspace{0.3cm}

These equations are immediate using $(6)$. Now we are ready to state the integration by parts.

\bt{Integration by parts for star-integrals}

\begin{equation}
    \starint f(x)^{g'(x)} \hspace{3pt} dx = f(x)^{g(x)}\starint f^{*}(x)^{-g(x)} \hspace{3pt} dx
\end{equation}

\et

\vspace{0.3cm}

which is in fact that same integral we had in the beginning of this section, now using the star-derivative. Also, we can discuss a ``Taylor series" for these functions. Consider

$$f(x) = \prod_{i=0}^{\infty}a_{i}^{(x-c)^i}$$

\vspace{0.3cm}

and we can try to compute $a_{i}$. Notice $f(c) = a_{0}$. Using part $(d)$ from remark 5.3, we have

$$f^{*}(x) = \prod_{i=1}^{\infty}a_{i}^{i(x-c)^{i-1}},$$

\vspace{0.3cm}

since the star-derivative of a constant is 1. Again, plugging in $x=c$, we have $f^{*}(c) = a_{1}$. Continuing this, we get the formula $\displaystyle a_{i} = \Big(f^{*(i)}(c)\Big)^{1/i!}$. And thus,

\vspace{0.2cm}

\bt[Taylor's theorem]{Given $f(x)$ is infinitely star-differentiable at $x=c$, we may write}

\begin{equation}
    f(x) = \prod_{i=0}^{\infty} \Big(f^{*(i)}(c)\Big)^{(x-c)^i/i!}
\end{equation}

where $\displaystyle f^{(i)*}(c) = \frac{\delta^i}{\delta x^{i}}(f(x))\Big|_{x=c}$ is the $i$-th star-derivative of $f(x)$ evaluated at $x=c$.

\et

\vspace{0.3cm}

One can discuss radius of convergence and intervals of convergence and more with this as well. Taking the logarithm of both sides, we transform the product into a sum and using the standard Taylor's theorem, we have

$$\frac{d^i}{dx^i}\Big(\log(f(x))\Big) = \log\Big(\frac{\delta^i}{\delta x^{i}}(f(x))\Big), \hspace{5pt} \forall \hspace{3pt} i \in \NN.$$

\vspace{0.3cm}

\section{Geometric interpretation}

\vspace{0.3cm}

We know that the regular derivative of a function represents the slope of the graph and the integral represents the area underneath the graph. Can we give some geometric interpretation to these star-integrals and star-derivatives? First, we provide the fundamental theorem of calculus for star-integrals.

\vspace{0.3cm}

\bt[Fundamental theorem of calculus for star-integrals]{Define}

$$F(x) := \starint f(x) \hspace{3pt} dx$$

then $F^{*}(x) = f(x).$ Further,

\begin{equation}
\starint_{\hspace{-5pt}a}^{b} f(x) \hspace{3pt} dx = \frac{F(b)}{F(a)}
\end{equation}

\et

\vspace{0.3cm}

If we try to relate this star-integral to the area under a graph, one method is to see if there exists a $g(x)$ such that $\displaystyle \starint f(x) \hspace{3pt} dx = \int g(x) \hspace{3pt} dx$. Naturally, $g(x)$ is the standard derivative of $\displaystyle \starint f(x) \hspace{3pt} dx$. We will have

$$g(x) = \frac{d}{dx}\bigg(\starint f(x) \hspace{3pt} dx\bigg) = \frac{d}{dx}\Bigg(\exp\bigg(\int \log(f(x)) \hspace{3pt} dx\bigg)\Bigg) = \log(f(x))\exp\bigg(\int \log(f(x))\bigg)$$

\vspace{0.3cm}

Let's try this out. Let $f(x)=x$ and our $g(x)$ quickly becomes $\displaystyle g(x) = \log(x)\Big(\frac{x}{e}\Big)^x$. However, if we integrate this from $0$ to $1$ for example,

$$\frac{1}{e} = \starint_{\hspace{-5pt}0}^{1} x \hspace{3pt} dx = \int_{0}^{1} \log(x)\Big(\frac{x}{e}\Big)^x \hspace{3pt} dx = \Big(\frac{x}{e}\Big)^x\bigg|_{0}^{1} = \frac{1}{e}-1$$

\vspace{0.3cm}

What happened? The problem is from $(9)$. The fundamental theorem of calculus for star-integrals divides the two antiderivatives by each other instead of subtracting. Instead we need to divide $1/e$ and $1$. But using regular integration, this is impossible. From this, it strongly hints that this star-integral cannot be made geometric, at least in the Euclidean sense. We would need to define a multiplicative metric and work in that multiplicative metric space, using non-Euclidean geometry. A multiplicative metric $d^{*}(x,y)$ is one such that 

\vspace{0.3cm}

$\displaystyle (a) \hspace{7pt} d^{*}(x,y) \geq 1$ and $d^{*}(x,y)=1 \Longleftrightarrow x=y$

\vspace{0.1cm}

$\displaystyle (b) \hspace{7pt} d^{*}(x,y) = d^{*}(y,x)$

\vspace{0.1cm}

$\displaystyle (c) \hspace{7pt} d^{*}(x,z) \leq d^{*}(x,y)d^{*}(y,z)$

\vspace{0.3cm}

for $x,y,z$ in the given space. These types of metrics have been studied (see \cite{Dot}) but there has been no connection between them and star-integrals or star-derivatives. Even though there does not seem to be anything geometric about these new operators, we can still discuss geometric theorems such as mean value theorem.

\bt[Mean value theorem for star-integrals]{If a function $f(x)$ is continuous on $[a,b]$ then $\exists \hspace{4pt} c \in (a,b)$ such that}

$$\starint_{\hspace{-5pt}a}^{b} f(x) \hspace{3pt} dx = f(c)^{b-a}.$$

\et

\vspace{0.3cm}

For example, $\displaystyle \starint_{\hspace{-5pt}1}^{2} e^{1/x} \hspace{3pt} dx = 2 = e^{1/c}$ when $c=1/\log(2) \in (1,2)$. Then for derivatives,

\bt[Mean value theorem for star-derivatives]{If a function $f(x) > 0$ is continuous on $[a,b]$ and differentiable on $(a,b)$ then $\exists \hspace{4pt} c \in (a,b)$ such that}

$$f^{*}(c) = \bigg(\frac{f(b)}{f(a)}\bigg)^{1/(b-a)}.$$

\et

\vspace{0.3cm}

Taking the logarithm of either of these equations, one arrives back at the regular mean value theorem statements. Also we can get Rolle's theorem too that if $f(a)=f(b)$, then there is a $c \in (a,b)$ such that $f^{*}(c)=1$.

\vspace{0.3cm}

\section{Extensions and generalizations}

We can keep discussing more things about these star-integrals. One idea could be ways to extend the star integral to negative-valued functions. Another option could be talking about Lebesgue star-integrable functions and perhaps discussing a norm $||.||_{p}^{*}$ and a space $L_{*}^{p}$, although it may need a different name so as to not confuse it with the dual space of $L^p$. But even before this, we may need to talk about improper star-integrals. It is interesting that from 3.1.d with $f(x) = 1$, we see $\displaystyle \starint_{\hspace{-5pt}0}^{x} k \hspace{3pt} dx = k^x$. So if we take a limit as $x\rightarrow\infty$, this star-integral converges for $0<k\leq1$ and diverges for $k>1$, which is quite strange. Even more strange is because $\log(f(x)) < f(x)$, if the regular integral is finite, it necessarily implies that the star-integral is finite. Perhaps one should define star-integrals that don't converge as those that approach infinity or 0.

\vspace{0.3cm}

We can also think of generalizations for this operator. The exponential and natural logarithm functions are inverses of each other. So we could define a different integral in a similar fashion but for different inverse functions. Let us define $\displaystyle \mathcal{L}(f(x),G(x)):=G\Bigg(\int G^{-1}(f(x)) \hspace{3pt} dx\Bigg) $. Then,

\vspace{0.3cm}

$$\starint f(x) \hspace{3pt} dx = \mathcal{L}(f(x),e^x)$$

and

$$\int f(x) \hspace{3pt} dx = \mathcal{L}(f(x),x).$$

\vspace{0.3cm}

Moving along, we can define

$$\int_{\#} f(x) \hspace{3pt} dx = \mathcal{L}(f(x),\log(x))$$

or

$$\int_{\$} f(x) \hspace{3pt} dx = \mathcal{L}(f(x),x^2)$$

where we take the positive square root. Perhaps these integrals also have some nice properties to them as well. A last generalization we can discuss is double star-integrals, for example

$$\starint \starint_{\hspace{-5pt}0}^{y} k \hspace{3pt} dx \hspace{2pt} dy = Ck^{y^2/2}$$

or

$$\starint \starint_{\hspace{-5pt}1}^{y} e^{k/x} \hspace{3pt} dx \hspace{2pt} dy = C\Big(\frac{y}{e}\Big)^{ky}$$

\vspace{0.3cm}

and we can perhaps find a similar Stokes' Theorem or Divergence Theorem. There are many different paths to take with star-integrals and star-derivatives and one may discover even more interesting properties.

\vspace{1cm}

\textbf{Acknowledgements.} I would like to thank Cristhian Emmanuel Garay L\'opez, Andre Issa, and Hetansh Mehta for their inspiring ideas and comments which served as motivation for this article.


\bigskip

\begin{flushright}
\begin{minipage}{148mm}\sc\footnotesize
University of Pittsburgh, Department of Mathematics, 301 Thackeray Hall, Pittsburgh, PA 15260, USA\\
{\it E--mail address}: {\tt djo15@pitt.edu} \vspace*{3mm}
\end{minipage}
\end{flushright}

\end{document}